\documentclass{amsart}

\usepackage{amsmath}
\usepackage{amsfonts}
\usepackage{amsthm}
\usepackage{setspace}
\usepackage{fancyhdr}
\usepackage{graphicx}
\usepackage{verbatim}
\usepackage[pdftex,colorlinks=true, pdfstartview=FitV, linkcolor=blue, citecolor=blue, urlcolor=blue]{hyperref}
\usepackage[T1]{fontenc}

\newtheorem{theorem}{Theorem}

\newtheorem{lemma}[theorem]{Lemma}

\newcommand{\Mat}{\ensuremath{\textrm{Mat}}}
\newcommand{\M}{\textit{Macaulay2}\ }

\newcommand{\defn}[1]{\textsf{#1}}
\newcommand{\QS}{\textit{QuillenSuslin}\ }

\title[QuillenSuslin]{The QuillenSuslin Package for Macaulay2}
\author[B. Barwick and B. Stone]{Brett Barwick and Branden Stone}
\date{\today}

\address{Brett Barwick, Division of Mathematics and Computer Science, University of South Carolina Upstate, Spartanburg, SC 29303}
\email{bbarwick@uscupstate.edu}
\thanks{The first author was partially supported by NSA Grant No. H98230-10-1-0361.}

\address{Branden Stone, Mathematics Program, Bard College/Bard Prison Initiative, P.O. Box 5000, Annandale-on-Hudson, NY 12504}
\email{bstone@bard.edu}

\subjclass[2010]{Primary: 13P99; Secondary: 13C10}

\begin{document}

\maketitle

\begin{abstract}
	The \QS package for \M provides the ability to compute a free basis for a projective module over a polynomial ring with coefficients in $\mathbb{Q}, \mathbb{Z},$ or $\mathbb{Z}/p$ for a prime integer $p$. A brief description of the underlying algorithm and the related tools are given.
\end{abstract}

\section{Introduction} 

	In 1955, J-P. Serre posed the following question: do there exist finitely generated projective $k[x_1,\ldots,x_n]$ modules, with $k$ a field, which are not free?~\cite{Serre}  This question was known as ``Serre's Problem'' and the question in its full generality remained open for 21 years until it was resolved independently by D. Quillen and A. A. Suslin in 1976, resulting in the following well-known theorem.

\begin{theorem}[Quillen-Suslin, 1976 \cite{Qui76,Sus76}]
Let $S = k[x_1,\ldots,x_n]$, with $k$ a field.  Then every finitely generated projective $S$-module is free.
\end{theorem}

	However, the proofs given were not entirely constructive, and it was not until the early 1990's that papers such as \cite{FG90, LS92, LW00} began giving fully constructive versions of the proof.  In 1992, A. Logar and B. Sturmfels \cite{LS92} published the algorithmic proof of the Quillen-Suslin Theorem that forms the basis for the methods in \emph{QuillenSuslin}.  In their paper, Logar and Sturmfels describe, via the technique of completion of unimodular rows, how to construct a free generating set for a projective module over $\mathbb{C}[x_1,\ldots,x_n]$.  One can extend these constructive techniques to work over more general coefficient rings such as $\mathbb{Q}, \mathbb{Z},$ and $\mathbb{Z}/p\mathbb{Z}$, for $p$ a prime integer.  Descriptions of some of these more general techniques, along with applications to areas such as systems control theory, can be found in \cite[Ch.~2]{Fab09} and \cite[Ch.~3]{FQ07}, and the algorithms given in these papers were used to create a similar QuillenSuslin package \cite{QSMaple} for the computer algebra system \emph{Maple}.  We have implemented these algorithms, with some modifications, in our \emph{QuillenSuslin} package for \emph{Macaulay2} \cite{M2}.  In the next section we will give some preliminary definitions and results that reduce the statement of the Quillen-Suslin Theorem to a more concrete matrix theoretic problem concerning the completion of unimodular rows over polynomial rings to square invertible matrices. 

\section{Preliminaries}

	In this section, $R$ will denote a commutative ring and $M$ will denote a finitely generated $R$-module.  We say that $M$ is a \defn{projective} $R$-module if it is a direct summand of a free module.  Similarly, we define a slightly stronger notion by saying that $M$ is \defn{stably free} if there exists some $m \geq 0$ such that $M \oplus R^m$ is free.  A module $M$ is stably free if and only if it is isomorphic to the kernel of a surjective $R$-linear map $\phi : R^n \to R^m$ for some $m \leq n$.  Since $\phi$ surjects onto a free module, we know that this map splits and admits a right inverse $\psi:R^m \to R^n$ so that $\phi \psi = \textrm{id}_{R^m}$.  Therefore a matrix representing $\phi$ is right invertible, and we call such a right invertible matrix over $R$ \defn{unimodular}.  Using this terminology, it is not difficult to show that the Quillen-Suslin Theorem as stated above is equivalent to the following matrix theoretic statement about unimodular matrices.

\begin{theorem}[Quillen-Suslin, restatement {\cite[Theorem 1.1]{LS92}}]\label{thm:QSrestate}
Let $S=R[x_1,\ldots,x_n]$ with $R$ a principal ideal domain and let $U \in \Mat_{m \times n}(R)$ be a unimodular matrix over $S$ with $m \leq n$.  Then there exists a unimodular matrix $V \in \Mat_{n \times n}(S)$ such that 
\[
	UV = \left[ \begin{array}{cccc|ccc} 
					1 & 0 & \cdots & 0 & 0 & \cdots & 0 \\
					0 & 1 & \cdots & 0 & 0 & \cdots & 0 \\
					 \vdots & \vdots & \ddots & \vdots & \vdots & \ddots & \vdots \\
					0 & 0 & \cdots & 1 & 0 & \cdots & 0 
				\end{array} \right].
\]
\end{theorem}

	A matrix $V$ that satisfies the properties in the above theorem is said to \defn{solve the unimodular matrix problem for $U$}.  Given such a $V$, the first $m$ rows of the invertible matrix $V^{-1}$ are the same as the original matrix $U$.  Therefore we see that proving the Quillen-Suslin Theorem is equivalent to showing that any unimodular matrix can be completed to a square invertible matrix over the polynomial ring.  As it turns out, it suffices to show that the \defn{unimodular row problem} can be solved, that is, Theorem \ref{thm:QSrestate} holds for unimodular row vectors \cite{LS92}.  For more details concerning this equivalent formulation of the Quillen-Suslin Theorem, we refer the interested reader to the excellent book of Lam \cite{Lam06}.

\section{The Logar-Sturmfels Algorithm} 

	Before describing the general algorithm, we mention that \emph{QuillenSuslin} contains several shortcut methods as described in \cite[Sect. 2.2]{Fab09}. These methods allow us to quickly solve the unimodular row problem for a row satisfying certain properties; often allowing us to avoid the worst-case general algorithm.  These shortcut methods are automatically used as soon as they are applicable during the methods \verb|computeFreeBasis|, \texttt{completeMatrix}, \texttt{qsAlgorithm}, and \texttt{qsIsomorphism}. \smallskip

	The main idea behind the Logar-Sturmfels algorithm is to iteratively reduce the number of variables involved in a unimodular row $\pmb{f}$ one by one, eventually obtaining a unimodular row $\tilde{\pmb{f}}$ over the coefficient ring, which is a PID.  We can then use a simple algorithm based on the Smith normal form of $\tilde{\pmb{f}}$ to construct a final unimodular matrix $U$ so that $\tilde{\pmb{f}} U = [1 \ 0 \ \cdots \ 0]$.  Multiplying together all of the matrices used during the process, one can construct a unimodular matrix that solves the unimodular row problem for the original row $\pmb{f}$. \smallskip

	The process of eliminating a variable from a unimodular matrix is organized into three main steps: the \defn{normalization step}, the \defn{``local loop''}, and the \defn{patching step}.  Below is a brief description of each step, as well as a demonstration of the corresponding commands in the \emph{QuillenSuslin} package.  We will work over the polynomial ring $S = \mathbb{Z}[x,y]$ and consider the unimodular row $\pmb{f} = [x^2 \, , \, 2y+1 \, , \, x^5y^2+y]$ over $S$. \smallskip

\begin{verbatim}
     i1 : loadPackage "QuillenSuslin";
     i2 : S = ZZ[x,y];

     i3 : f = matrix {{x^2,2*y+1,x^5*y^2+y}}
     o3 = | x2 2y+1 x5y2+y |
                  1       3
     o3 : Matrix S  <--- S

\end{verbatim}

\indent We can use the command \texttt{isUnimodular} to check that this row is indeed unimodular over $S$.
\begin{verbatim}
     i4 : isUnimodular f
     o4 = true
\end{verbatim}

	In order to eliminate the variable $y$ from this unimodular row, we will construct a unimodular matrix $U$ so that $\pmb{f}U$ is the same as $\pmb{f}$ with $y$ replaced by 0.  We first demonstrate the normalization step.

\subsection{Normalization Step}

	Since Horrocks' Theorem (see Theorem \ref{thm:Horrocks} below) requires a monic polynomial, we must first construct a unimodular matrix $U_1$ and an invertible change of variables $x_i \leftrightarrow X_i$ so that the first entry of $\pmb{f}U_1$ is monic in $X_n$.  The normalization step is based on the following result, which has been slightly restated for our purposes.

\begin{lemma}[{\cite[Lemma 10.6]{SV76}}]
Let $R$ be a Noetherian ring, $n$ and $m$ natural numbers, $S = R[x_1,\ldots,x_n]$, $m \geq \dim R + 2$, and $\pmb{f} = [f_1 \ f_2 \ \cdots \ f_m]$ a unimodular row over $S$.  Then there exists an $m \times m$ unimodular matrix $U$ over $S$, and an invertible change of variables $x_i \leftrightarrow X_i$, so that after applying the change of variables the first entry of $\pmb{f} U$ is monic in $X_n$ when viewed as a polynomial in $R[X_1,\ldots,X_{n-1}][X_n]$.
\end{lemma}

	A constructive version of this result is implemented in the method \texttt{changeVar}, and is used as in the following example.
\begin{verbatim}
     i5 : (U1,subs,invSubs) = changeVar(f,{x,y})
     o5 = (| 1 0 0 |, | y x |, | y x |)
           | 0 1 0 |
           | 0 0 1 |
     o5 : Sequence	

     i6 : f = sub(f*U1,subs)
     o6 = | y2 2x+1 x2y5+x |
                  1       3
     o6 : Matrix S  <--- S
\end{verbatim}
Notice that since the first entry of the row $\pmb{f}$ was already monic in $x$, the method simply returned a permutation of the variables interchanging $x$ and $y$ so that the first entry of the new row would be monic in $y$.  Now that the first entry of the row is monic in the variable we are trying to eliminate, we may proceed to the ``local loop.''

\subsection{Local Loop}

	The purpose of the local loop is to compute a collection of \emph{local solutions} to the unimodular row problem for $\pmb{f}$.  The local loop is based on the following result of Horrocks.

\begin{theorem}[Horrocks, {\cite[Prop. 4.98]{Rot09}}]\label{thm:Horrocks}
Consider the polynomial ring $B[y]$, where $B$ is a local ring, and let $\pmb{f} = [f_1 \ f_2 \ \cdots \ f_m]$ be a unimodular row over $B[y]$.  If some $f_i$ is monic in $y$, then there exists a unimodular $m \times m$ matrix $U$ over $B[y]$ so that $\pmb{f} U = [1 \ 0 \ \cdots \ 0]$.
\end{theorem}

	In order to eliminate the last variable $x_n$ in a unimodular row over a polynomial ring $R[x_1,\ldots,x_n]$, the local loop proceeds in the following way:  First set $I=(0)$ in $A = R[x_1,\ldots,x_{n-1}]$.  Now while $I \not = A$, the $i^{\text{th}}$ iteration of the loop is
\begin{enumerate}
	\itemsep 1pt
	\parskip 0pt
	\item Find a maximal ideal $\mathfrak{m}_i$ in $A$ containing $I$.
	\item Apply Horrocks' Theorem to the row $\pmb{f}$, viewed as a unimodular row over $A_{\mathfrak{m}_i}[x_n]$, to find a unimodular matrix $L_i$ over $A_{\mathfrak{m}_i}[x_n]$ that solves the unimodular row problem for $\pmb{f}$ (we call this $L_i$ a \defn{local solution to the unimodular row problem for} $\pmb{f}$).
	\item Let $d_i$ denote the common denominator for all of the elements in the matrix $U_i$.
	\item Set $I = I + (d_i)$.
\end{enumerate}
If $I \not = A$, then we repeat the loop.  Otherwise we are able to stop and go on to the patching step.  Notice that since we are creating a strictly ascending chain of ideals $(d_1) \subset (d_1,d_2) \subset \cdots$ in the Noetherian ring $A$, this loop must terminate in a finite number of steps with $(d_1,\ldots,d_k) = A$ for some integer $k$. \smallskip

	In our example, we use the method \texttt{getMaxIdeal} to first find an arbitrary maximal ideal in $\mathbb{Z}[x]$, and we set $\mathfrak{m}_1 = (2,x)$.  Using the method \texttt{horrocks}, we can compute a unimodular matrix $L_1$ over $\left( \mathbb{Z}[x]_{(2,x)} \right)[y]$ so that $\pmb{f}L_1 = [1 \ 0 \ 0]$.
\begin{verbatim}
     i7 : A = ZZ[x];

     i8 : m1 = getMaxIdeal(ideal(0_A),{x}) 
     o8 = ideal (2, x)
     o8 : Ideal of A

     i9 : L1 = horrocks(f,y,m1)
     o9 = | 0        1          0                |
          | 1/(2x+1) -y2/(2x+1) (-x2y5-x)/(2x+1) |
          | 0        0          1                |
                         3              3
     o9 : Matrix (frac S)  <--- (frac S)
\end{verbatim}
Since $d_1 = 2x+1$ is a common denominator for the entries of $L_1$ and $(2x+1) \not = \mathbb{Z}[x]$, we use \texttt{getMaxIdeal} again to find a maximal ideal containing $2x+1$, and we set $\mathfrak{m}_2 = (3,x-1)$.  We use \texttt{horrocks} a second time to compute a new local solution \texttt{L2} with common denominator $d_2 = x$.
\begin{verbatim}
     i10 : m2 = getMaxIdeal(sub(ideal(2*x+1),A),{x})
     o10 = ideal (x - 1, 3)
     o10 : Ideal of A
     
     i11 : L2 = horrocks(f,y,m2) 
     o11 = | -xy3 xy5+1 2x2y3+xy3 |
           | 0    0     1         |
           | 1/x  -y2/x (-2x-1)/x |

                          3              3
     o11 : Matrix (frac S)  <--- (frac S)                                                                                                                                                       


     i12 : sub(ideal(2*x+1,x),S) == ideal(1_S)
     o12 = true
\end{verbatim}
Since $(d_1,d_2) = (2x+1,x) = \mathbb{Z}[x]$, we are able to exit the local loop and proceed to the patching step.

\subsection{Patching Step}	

	Loosely speaking, the patching step involves multiplying slight variations of the local solutions $L_1,\ldots,L_k$ together in a clever way so that the product $U$ is a unimodular matrix over the polynomial ring $R[x_1,\ldots,x_n]$ and multiplying $\pmb{f}$ times $U$ is equivalent to evaluating $\pmb{f}$ when $x_n = 0$, thereby eliminating one of the variables in the row $\pmb{f}$.  For more details, see \cite[pg.~235]{LS92}. \smallskip
	
	Following along with our example, we use the method \texttt{patch} applied to our list $\{L_1,L_2\}$ of local solutions and we specify that $y$ is the variable that we want to eliminate.
\begin{verbatim}
     i13 : U = patch({L1,L2},y)
     o13 = | -32x6y5+1                  0 8x5y3             |
           | 16x5y7-8x4y7+4x3y7+2xy2-y2 1 -4x4y5+2x3y5-x2y5 |
           | -4xy2                      0 1                 |
                   3       3
     o13 : Matrix S  <--- S                                                                                                                                                                                          

     i14 : f*U
     o14 = | 0 2x+1 x |
                   1       3
     o14 : Matrix S  <--- S                                                                                                                                                                                          
\end{verbatim}
We can see that multiplying the row $\pmb{f}$ times the unimodular matrix \texttt{U} is equivalent to evaluating $\pmb{f}$ when $y=0$ (keeping in mind that the variables $x$ and $y$ were interchanged during the normalization step).

\section{Core methods in the QuillenSuslin package} 

The method \texttt{qsAlgorithm} automates all of the above computations for computing a solution to the unimodular matrix problem, and automatically applies the shortcut methods in \cite{Fab09} when possible.  We demonstrate the use of \texttt{qsAlgorithm} by finding a solution to the unimodular row problem for the row $\pmb{f} = [y^2 \, , \, 2x+1 \, , \, x^2y^5+x]$ given earlier.
\begin{verbatim}
     i15 : U = qsAlgorithm f
     o15 = | 2x3y2 -2x5y2+1 -2x8y4-2x3y3 |
           | 1     -x2      -x5y2-y      |
           | -2    2x2      2x5y2+2y+1   |
                   3       3
     o15 : Matrix S  <--- S                                                                                                                                                                                          

     i16 : det U
     o16 = -1
     o16 = S
     i17 : f*U
     o17 = | 1 0 0 |
                   1       3
     o17 : Matrix S  <--- S                                                                                                                                                                                          
\end{verbatim}

	The package also contains a method \texttt{completeMatrix} that completes a unimodular matrix over a polynomial ring to a square invertible matrix.  Again we demonstrate its use on the unimodular row $\pmb{f}$.
\begin{verbatim}
     i18 : C = completeMatrix f
     o18 = {0} | y2 2x+1   x2y5+x |
           {2} | 1  -2x2y3 0      |
           {7} | 0  2      1      |
                   3       3
     o18 : Matrix S  <--- S                                                                                                                                                                                          

     i19 : det C
     o19 = -1
     o19 : S
\end{verbatim}

	The previous two methods, \texttt{qsAlgorithm} and \texttt{completeMatrix}, also work over Laurent polynomial rings of the form $k[x_1^{\pm 1},\ldots,x_n^{\pm 1}]$ with $k = \mathbb{Q}$ or $\mathbb{Z}/p\mathbb{Z}$ for $p$ a prime integer.  The algorithm in the Laurent polynomial case is due to Park and is described in \cite[p.~215]{Park}.

        Finally, we give an example to demonstrate the method \texttt{computeFreeBasis} that computes a free generating set for a projective module.  We define $K = \ker \pmb{f}$, which we can check is a projective $\mathbb{Z}[x,y]$-module by using the command \texttt{isProjective}.
\begin{verbatim}
     i20 : K = ker f
     o20 = image {2} | 2y+1 2x3y3+x3y2 2x5y2-1 -x5y2-y |
                 {1} | -x2  y          x2      0       |
                 {7} | 0    -2y-1      -2x2    x2      |
                                   3
     o20 : R-module, submodule of R

     i21 : isProjective K
     o21 = true

     i22 : mingens K
     o22 = {2} | -2y-1 2x3y3+x3y2 x5y2+y |
           {1} | x2    y          0      |
           {7} | 0     -2y-1      -x2    |
                   3       3
     o22 : Matrix R  <--- R

     i23 : syz mingens K
     o23 = {3} | y    |
           {8} | -x2  |
           {9} | 2y+1 |
                   3       1
     o23 : Matrix R  <--- R

     i24 : phi = qsIsomorphism K
     o24 = {3} | 0 0 |
           {8} | 1 0 |
           {9} | 0 1 |
           {9} | 0 0 |
     o24 : Matrix

     i25 : source phi
             2
     o25 = R
     o25 : R-module, free

     i26 : target phi
     o26 = image {2} | 2y+1 2x3y3+x3y2 2x5y2-1 -x5y2-y |
                 {1} | -x2  y          x2      0       |
                 {7} | 0    -2y-1      -2x2    x2      |
                                   3
     o26 : R-module, submodule of R

     i27 : isIsomorphism phi
     o27 = true

     i28 : B = computeFreeBasis K
     o28 = {2} | 2x3y3+x3y2 2x5y2-1 |
           {1} | y          x2      |
           {7} | -2y-1      -2x2    |
                   3       2
     o28 : Matrix R  <--- R

     i29 : syz B
     o29 = 0
                   2
     o29 : Matrix R  <--- 0

     i30 : image B == K
     o30 = true
\end{verbatim}
From the \emph{Macaulay2} output, we can see that the native command \texttt{mingens} does not produce a free generating set for $K$, while \texttt{computeFreeBasis} produces a set of 2 generators for $K$ with no relations, demonstrating that $K$ is free.

\section{Acknowledgements} 

	We thank Amelia Taylor for organizing the 2010 \textit{Macaulay2} workshop at Colorado College\footnote{The \M workshop at Colorado College in Colorado Springs, Colorado was funded by the the National Science Foundation under Grant No. 09-64128 and the National Security Agency under Grant No. H98230-10-1-0218.} where work on the \QS package began.  We would especially like to thank Hirotachi Abo for many useful conversations during the workshop as well as Jason McCullough.  
	
	The first author thanks his advisor, Andrew R. Kustin, for many useful discussions during the development of the package.


\begin{thebibliography}{Lam06}

\bibitem[Fab]{QSMaple}
  Anna Fabianska, \emph{QuillenSuslin Package for Maple}, Available at \href{http://wwwb.math.rwth-aachen.de/QuillenSuslin}%
  {http://wwwb.math.rwth-aachen.de/QuillenSuslin}.

\bibitem[Fab09]{Fab09}
Anna Fabianska, \emph{Algorithmic analysis of presentations of groups and
  modules}, Ph.D. thesis, RWTH Aachen University, Jan 2009.

\bibitem[FQ07]{FQ07}
Anna Fabianska and Alban Quadrat, \emph{Applications of the {Q}uillen-{S}uslin Theorem to multidimensional systems theory}, Gro\"obner bases in control theory and signal processing, Radon Ser. Comput. Appl. Math. \textbf{3}, 23--106. \MR{2402709 (2009f:13013)}

\bibitem[FG90]{FG90}
Noa{\"{\i}} Fitchas and Andr{\'e} Galligo, \emph{Nullstellensatz effectif et
  conjecture de {S}erre (th\'eor\`eme de {Q}uillen-{S}uslin) pour le calcul
  formel}, Math. Nachr. \textbf{149} (1990), 231--253. \MR{1124807 (92i:12002)}

\bibitem[GS]{M2}
Daniel~R. Grayson and Michael~E. Stillman, \emph{Macaulay2, a software system
  for research in algebraic geometry}, Available at
  \href{http://www.math.uiuc.edu/Macaulay2/}%
  {http://www.math.uiuc.edu/Macaulay2/}.

\bibitem[Lam06]{Lam06}
T.~Y. Lam, \emph{Serre's problem on projective modules}, Springer Monographs in
  Mathematics, Springer-Verlag, Berlin, 2006. \MR{2235330 (2007b:13014)}

\bibitem[LW00]{LW00}
Reinhard~C. Laubenbacher and Cynthia~J. Woodburn, \emph{A new algorithm for the
  {Q}uillen-{S}uslin theorem}, Beitr\"age Algebra Geom. \textbf{41} (2000),
  no.~1, 23--31. \MR{1745576 (2001e:13013)}

\bibitem[LS92]{LS92}
Alessandro Logar and Bernd Sturmfels, \emph{Algorithms for the
  {Q}uillen-{S}uslin theorem}, J. Algebra \textbf{145} (1992), no.~1, 231--239.
  \MR{1144671 (92k:13006)}

\bibitem[Par04]{Park}
Hyungju Park, Symbolic computation and signal processing, J. Symbolic Comput. 37
(2004), no.~2, 209--226. \MR{2093232 (2005d:94039)}

\bibitem[Qui76]{Qui76}
Daniel Quillen, \emph{Projective modules over polynomial rings}, Invent. Math.
  \textbf{36} (1976), 167--171. \MR{0427303 (55 \#337)}

\bibitem[Rot09]{Rot09}
Joseph~J. Rotman, \emph{An introduction to homological algebra}, second ed.,
  Universitext, Springer, New York, 2009. \MR{2455920 (2009i:18011)}

\bibitem[Ser55]{Serre}
Jean-Pierre Serre, \emph{Faisceaux alg\'ebriques coh\'erents}, Ann. of Math.
  (2) \textbf{61} (1955), 197--278. \MR{0068874 (16,953c)}

\bibitem[Sus76]{Sus76}
A.~A. Suslin, \emph{Projective modules over polynomial rings are free}, Dokl.
  Akad. Nauk SSSR \textbf{229} (1976), no.~5, 1063--1066. \MR{0469905 (57
  \#9685)}

\bibitem[VS76]{SV76}
LN~Vaserstein and A.~A. Suslin, \emph{{Serre's Problem on Projective Modules
  over Polynomial Rings, and Algebraic K-Theory}}, {Mathematics of the
  USSR-Izvestiya} \textbf{{10}} ({1976}), no.~{5}, {937--1001} ({English}).

\end{thebibliography}


\providecommand{\bysame}{\leavevmode\hbox to3em{\hrulefill}\thinspace}
\providecommand{\MR}{\relax\ifhmode\unskip\space\fi MR }
\providecommand{\MRhref}[2]{%
  \href{http://www.ams.org/mathscinet-getitem?mr=#1}{#2}
}
\providecommand{\href}[2]{#2}


\end{document}